\documentclass{amsart}

\title[Riemann Surfaces]{Embedding Riemann Surfaces Properly into $\CC^2$}
\author{Erlend Forn\ae ss Wold}
\date{September 22, 2005}

\usepackage{amscd}

\newtheorem{theorem}{Theorem}
\newtheorem{lemma}{Lemma}
\newtheorem{proposition}{Proposition}

\theoremstyle{definition}

\theoremstyle{remark}
\newtheorem{remark}{Remark}

\newcommand{\NN}{\mathbb{N}}
\newcommand{\RR}{\mathbb{R}}
\newcommand{\CC}{\mathbb{C}}

\def\a{{\alpha}}

\def\d{{\delta}}
\def\e{{\epsilon}}
\def\s{{\sigma}}

\def\Aut{{\mathrm{Aut}}}

\begin{document}

\bibliographystyle{plain}

\begin{abstract}
For certain bordered submanifolds $M\subset\CC^2$  we show that
$M$ can be embedded properly and holomorphically into $\CC^2$.  An
application is that any subset of a torus with two boundary
components can be embedded properly into $\CC^2$.
\end{abstract}

\maketitle

\section{Introduction, Main Results and Notation}

In this paper we consider the problem of embedding bordered
Riemann surfaces properly into $\CC^2$.  A bordered Riemann
surface is obtained by taking a compact Riemann surface $\mathcal
R$, and removing a finite set of disjoint closed connected
components $D_1,...,D_m$, i.e. the bordered surface is
$\tilde{\mathcal R} =\mathcal R \setminus(\cup_{i=1}^m D_i)$.  The
case $\mathcal R=\CC\cup\{\infty\}$  and $m$  finite, was settled in
\cite{wd2}. The only case known thus far when the genus of
$\mathcal R$  is greater than zero, is the case where
$\tilde{\mathcal R}$ is a surface that is hyperelliptic (meaning
that its double is hyperelliptic), due to \v{C}erne and
Forstneri\v{c} \cite{cf}.  An example is the torus take away one
disc. This was proved by adjusting an embedding of
$\tilde{\mathcal R}$ into the polydisc (see Gouma \cite{go}
 and Rudin \cite{ru}), and composing with a certain
Fatou-Bieberbach map (as constructed by Stens\o nes \cite{St}  and
Globevnik \cite{Gl}). This naturally suggests the following
splitting of the general embedding problem into the following two
problems:

\

(a) Find some conditions on a bordered submanifold M of $\CC^2$
that enables you to embed it properly into $\CC^2$. \

(b) Embed $\tilde{\mathcal R}$  onto a surface $M\subset\CC^2$
satisfying the conditions from (a).

\

In \cite{cf} the authors established a condition that $M$
 is a certain kind of closed submanifold of the polydisc as sufficient, and they also
 suggested polynomial convexity as a condition on $M$.
\

In proving Theorem \ref{main} we establish another condition on
the submanifold $M$ (it follows from Proposition 3.1 in \cite{cf}
that a surface satisfying the condition established in \cite{cf}
can be perturbed to satisfy our condition), and in proving Theorem
\ref{torus} we solve the problem (b) for subsets of the torus with
two boundary components. Such a subset is not typically
hyperelliptic, so the result does not follow from \cite{cf}. In
proving Theorem \ref{jordan}, we show that by a linear change of
coordinates, the condition is always satisfied for a bounded
submanifold of $\CC^2$ whose boundary is a smooth Jordan curve
(such a surface is always polynomially convex).  In particular,
this implies the embedding theorem for any subset of a surface
with one boundary component whose double is hyperelliptic (via the
embedding provided by Rudin). Lastly we show that any bounded
submanifold of $\CC^2$ that is Runge can be exhausted by surfaces
admitting proper embeddings.

\begin{theorem}\label{main}
Let $M\subset\CC^2$  be a Riemann surface whose boundary
components are smooth Jordan curves $\partial_1,...,\partial_m$.
Assume that there are points $p_i\in\partial_i$  such that
$$
\pi_1^{-1}(\pi_1(p_i))\cap\overline M=p_i.
$$
Assume that $\overline M$  is a smoothly embedded surface, and that all $p_{i}$  are
regular points of the projection $\pi_{1}$.  Then $M$ can be properly holomorphically embedded
into $\CC^2$.
\end{theorem}
\begin{theorem}\label{torus}
Any subset of a torus with two boundary components can be properly
holomorphically embedded into $\CC^2$.
\end{theorem}
We note that one or both boundary components may reduce to a point.  The precise
regularity of the boundary of $D_i$ is not so important, as it is
always possible to embed $\tilde T$ onto a subset
of a torus such that the boundary consists of two real analytic components.
\begin{theorem}\label{jordan}
Let $M\subset\CC^2$  be a bounded Riemann surface such that
$\partial M$  is a smooth Jordan curve.  Then $M$ can be properly holomorphically
embedded into $\CC^2$.
\end{theorem}
In the following theorem we say that a surface $M\subset\CC^2$  is
Runge if for any compact set $K\subset\subset M$,  we have that
$\widehat K\cap\overline M\subset\subset M$, where $\widehat K$  denotes the
polynomially convex hull of $K$.
\begin{theorem}\label{exh}
Let $M\subset\subset\CC^2$  be a Riemann surface whose boundary is
a finite collection of smooth Jordan curves, and assume that $M$
is Runge.  Then $M$  has an exhaustion $M_i$  of Riemann surfaces
that are diffeomorphic to $M$, and each $M_i$  embeds properly holomorphically
into $\CC^2$.
\end{theorem}

\

As usual we will let $\mathbb{B}$  and $\triangle$  denote the
unit ball in $\CC^2$  and the unit disk in $\CC$  respectively,
and we let $B_R$  ($\triangle_R$) denote the ball (disk) centered
at the origin with radius $R$. We let $\pi_i$ denote the
projection on the ith coordinate axis, and we let $\|\cdot\|_K$
denote the euclidian sup-norm over a compact set $K$.  We let
$Aut_0(\CC^2)$ denote the group of holomorphic automorphisms
fixing the origin. Recall the definition of a basin of attraction:
If $\{F_j\}\subset Aut_0(\CC^2)$  is a countable sequence of
automorphisms fixing the origin, we will let $F(j)$ denote the
composition map $F_j\circ F_{j-1}\circ\cdot\cdot\cdot F_1$, and we
define the basin of attraction:
$$
\Omega_{\{F_j\}}^0=\{x\in\CC^2;\lim_{j\rightarrow\infty}
F(j)(x)=0\}
$$
A Fatou-Bieberbach domain is a proper sub-domain
$\Omega\subset\CC^2$, with a biholomorphic map
$\psi:\Omega\rightarrow\CC^2$ that is onto.

\

\verb"Acknowledgement:" The author would like to thank Franc
Forstneri\v{c} for directing his attention to the paper \cite{cf}.

\

\section{Fatou-Bieberbach domains and Proper Holomorphic Embeddings}

\

Let us start with a surface $M\subset\CC^2$  for which we want to
construct a holomorphic map
$$
\psi:M\rightarrow\CC^2
$$
that embeds $M$  properly into $\CC^2$.  It is sometimes
convenient to achieve this by constructing a Fatou-Bieberbach
domain $\psi\colon\Omega\rightarrow\CC^2$.  It is clear that
$\psi$ embeds $\Omega$  properly into $\CC^2$, so if
$M\subset\Omega$ and if $\partial M\subset\partial\Omega$, then
$\psi$  will embed $M$ properly.  Let us generalize this
situation.\

Assume that we are in the following setting:  We have two disjoint
sets $V,M\subset\CC^k$  and we want to construct a
Fatou-Bieberbach domain $\Omega$  such that
$M\subset\Omega\subset\CC^k\setminus V$. This is of course not
always possible, but we will give some conditions that are
sufficient for the construction to be possible. The conditions
will of course be such that we may apply them to prove Theorem
\ref{main}. \

\emph{Condition 1}: Let $K\subset\CC^k\setminus V$  be an
arbitrary polynomially convex compact set.  For any $R\in\RR^+$
and any $\e>0$ there exists an automorphism $\phi\in \Aut(\CC^k)$
such that the following are satisfied: \

(i) $\|\phi(x)-x\|<\e$  for all $x\in K$, \

(ii) $\phi(V)\subset\CC^k\setminus B_R$. \

\emph{Condition 2}: The set $M$ can be written as an increasing
union of compact sets $M=\cup_{i=1}^\infty K_i$  such that if
$K\subset\CC^k\setminus V$  is an arbitrary polynomially convex
compact set, then
$$
\widehat{K\cup K_j}\cap V=\emptyset
$$
for all $j\in\NN$. \

We have the following theorem:

\begin{theorem}\label{fb} Let $M,V\subset\CC^k$  be disjoint sets
satisfying the two conditions above.  Then there exists a
Fatou-Bieberbach domain $\Omega$  such that
$$
M\subset\Omega\subset\CC^k\setminus V
$$
\end{theorem}
\begin{proof}
We may assume that $\overline{\mathbb{B}}\cap V=\emptyset$.  Let
$A\colon\CC^k\rightarrow\CC^k$  be the linear automorphism defined
by
$$
A(z_1,...,z_k)=(\frac{z_1}{2},...,\frac{z_k}{2}).
$$
By Theorem 4 in \cite{wd1}  there exists a $\d>0$  such that if
$\{\s_j\}_{j=1}^\infty\subset \Aut_0(\CC^k)$  is a sequence of
automorphisms satisfying
$$
(*) \ \|\s_i(x)-A(x)\|<\d \ \mathrm{for \ all} \ x\in\mathbb{B},
$$
then the basin of attraction of the sequence $\s_j$ is
biholomorphic to $\CC^k$. \

We will prove the result by an inductive argument, and we let the
following be our induction hypothesis $I_j$:  We have a collection
of automorphisms $\{F_i\}_{i=1}^j\subset \Aut_0(\CC^k)$  with the
following satisfied: \

(a) Each $F_i$ is a finite composition of automorphisms $\s_k$
satisfying $(*)$, \

(b) $F(j)(K_j)\subset\mathbb{B}$,\

(c) $F(j)(V)\cap\overline{\mathbb{B}}=\emptyset$. \

We may assume that $I_1$  is satisfied with $F_1=A$.  By Condition
2 it follows that for $K=\overline{\mathbb{B}}\cup F(j)(K_{j+1})$
we have that
$$
\widehat K\cap F(j)(V)=\emptyset
$$
Let $s\in\NN$  such that $A^s(K)\subset\mathbb{B}$. For any
$\e>0$, by Condition 1 there exists an automorphism $\varphi\in
\Aut_0(\CC^k)$ such that \

(i) $\|\varphi(x)-x\|<\e$  for all $x\in K$, \

(ii) $\varphi(F(j)(V))\cap A^{-s}(B_2)=\emptyset$. \

If $\e$  is chosen small enough we have that
$$
\|A\circ\varphi(x)-A(x)\|<\d \ \mathrm{for \ all} \
x\in\overline{\mathbb{B}},
$$
so the composition map $F_{j+1}=A^s\circ\varphi$  clearly
satisfies (a).  It is easy to see that (b) and (c) are also
satisfied, so $F_{j+1}$  gives us $I_{j+1}$.  It follows that
$\Omega_{\{F_j\}}^0$  satisfies the claims of the theorem.
\end{proof}

The reason for proving this theorem is that Condition 1 and
Condition 2 can be proved to be satisfied for certain Riemann
surfaces $M\subset\CC^2$, and their boundaries $V=\partial M$.
Establishing this for a set of surfaces containing (up to
biholomorphism) all finitely connected subsets of $\CC$ was indeed
the content of the paper \cite{wd2}. \

In what follows, the following lemma from \cite{wd2}  will be an
essential ingredient in establishing Condition 1:
\begin{lemma}\label{move} Let $K\subset\CC^2$  be a polynomially
convex compact set, let $\e>0$, and let
$\Gamma=\{\gamma_j(t);j=1,..m, t\in [0,\infty)\}$ be a collection
of disjoint smooth curves in $\CC^2\setminus K$ without
self-intersection, such that
$lim_{t\rightarrow\infty}|\pi_1(\gamma_j(t))|=\infty$ for all $j$.
Assume that there exists an $N\in\RR$  such that
$\CC\setminus(\overline\triangle_R\cup\pi_1(\Gamma))$ does not
contain any relatively compact components for $R\geq N$.  Let
$p\in K$. Then for any $R\in\RR$  there exists an automorphism
$\phi\in Aut(\CC^2)$ such that the following is satisfied:

\

(i) $\|\phi(x)-x\|<\e$ for all $x\in K$, \

(ii) $\phi(\Gamma)\subset\CC^2\setminus B_R$,\

(iii) $\phi(p)=p$.
\end{lemma}
\section{Proof of Theorem \ref{main}}

\begin{proposition}\label{ext}
Let $M$  be a bordered submanifold of $\CC^2$  whose boundary is a
set of smooth curves $\partial_1,...,\partial_m$  that are all
unbounded. Then Condition 2 is satisfied for the pair $M$ and
$V=\partial M$.
\end{proposition}
\begin{proof}
Let $\{K_j\}$  be an exhaustion of $M$  by compact sets such that
the boundary of each $K_i$  is a finite collection of smooth
Jordan curves, and such that each $K_i$  is holomorphically convex
relative to $M$.  Choose an $m\in\NN$ such that $K\cap M\subset
K_i$ for all $i\geq m$. We claim that
$$
\widehat{K\cup K_i}=\widehat{K\cup\partial K_i}=K\cup K_i
$$
for all $i\geq m$.   By \cite{sb}  we have that $\widehat{K\cup
\partial K_i}\setminus(K\cup\partial K_i)$  is a closed subvariety of
$\CC^2\setminus(K\cup\partial K_i)$.  Observe that this  variety
does not contain isolated points. \

Assume that this set contains a variety $X$ different from
$K_i\setminus(K\cup\partial K_i)$, such that $\widehat{K\cup
K_i}=K\cup K_i\cup X$. We will show that this implies
$X=\emptyset$. \

We claim first that $\overline X\cap(K_i\setminus K)$ is at most a
discrete set of points in $K_i\setminus K$.  Let $x_0\in\partial
K_i\cap\overline X$, and observe that $Y:=\widehat{K_{i+1}\cup K}$
is a variety at $x_0$ (we may assume that $K_i\subset\subset
K_{i+1}$), and that $X\subset Y$.  Choose a small ball $B_\e(x_0)$
such that $Y_{x_0}:=Y\cap B_\e(x_0)$  may be written
$$
Y_{x_0}= S_0\cup S_1\cup\cdot\cdot\cdot\cup S_m;
$$
a union of irreducible components such that $Y_{x_0}$  is a
submanifold of $B_\e(x_0)$  except possibly at $x_0$. We arrange this
so that $S_0=B_\e(x_0)\cap K_{i+1}$.  We have that $X\cap
S_0=\emptyset$, for if not we would have $X\cap S_0=S_0\setminus
K_i$, which would imply that $X$  is unbounded, by the hypothesis
on $M$  and the fact that $K_{i}$  is holomorphically convex.
  This follows from Proposition 1, page 61 in \cite{ch}.
Now we have that
$$
\overline X\cap\partial K_i\cap
B_\e(x_0)\subset\overline{\cup_{j=1}^m S_j}\cap\partial
K_i=\{x_0\},
$$
and it follows that $x_0$  is an isolated point in the
intersection.  For a point $x_0\in K_i\setminus K$  not in the
boundary, repeat the same argument, this time using the assumption that
$X\cap K_{i}=\emptyset$,  and it follows that $x_0$  is isolated.  We have thus shown
that $\overline X\cap(K_i\setminus K)$  at most is a set
$P:=\{p_j\}_{j=1}^\infty$  which is discrete in $K_i\setminus K$.
In other words; the only accumulation points for the set $P$ in
$\CC^2$  has to be contained in $K$. \

By the local maximum modulus principle we have that
$$
\widehat{K\cup K_i}=(K\cup K_i)\cup\widehat{(\overline X\cap(K\cup
K_i))},
$$
and we have just shown that
$$
\overline X\cap(K\cup K_i)\subset K\cup P.
$$
Now we claim that the set $K\cup P$  is polynomially convex. Since
$K$  is polynomially convex, it has a Runge and Stein neighborhood
basis $\{\Omega_j\}$.  For any $j$  there is an $N\in\NN$  such
that $P_N:=\cup_{i=N}^\infty\{p_i\}\subset\Omega_j$.  Since
$\widehat{K\cup P_N}\subset\Omega_j$, and since the union of a
polynomially convex set and a finite set of points is polynomially
convex, it follows that $\widehat{K\cup P}\subset\Omega_j\cup P$.
Now let $j\rightarrow\infty$, and  we have shown that
$$ \widehat{(\overline X\cap(K\cup K_i))}\subset\widehat{K\cup
P}=K\cup P\Rightarrow\widehat{K\cup K_i}=K\cup K_i \ (\Rightarrow
X=\emptyset).
$$

\end{proof}

\begin{proposition}\label{em}
Let $M$  be a surface as in Theorem \ref{main}.  There exists an
embedding
$$
\varphi\colon M\rightarrow\CC^2
$$
such that Conditions 1 and 2 are satisfied for the pair
$\varphi(M)$ and $\partial(\varphi(M))$.
\end{proposition}
\begin{proof}
Let $c_1,...,c_m$  denote the derivatives $\frac{\partial
z_i}{\partial t}(0)$, and define
$$
\varphi(z,w)=(z,w+\sum_{i=1}^m\frac{a_i}{z-z_i(0)}).
$$

We will show that we may choose the coefficients $a_i$  such that
Condition 1 follows from Lemma \ref{move} (we look at $\pi_2$
instead of $\pi_1$). Let $\Gamma_i$ denote
$\varphi(\partial_i\setminus p_i)$ close to $p_i$, and let
$\Gamma_i^+$  and $\Gamma_i^-$ denote the parts where $t$  is
positive and negative.  As $t$ approaches $0$, we have that
$\pi_2(\Gamma_i^-)$ stays close to the line
$$
w(0)+\frac{a_i}{c_i t}.
$$
We also have that the absolute value of $\pi_2(\Gamma_i^-)$  is
strictly increasing when $t$ is close to and approaches zero. The
same conclusion holds for $\Gamma_i^+$ - and we observe that the
projection of the curve will point in the opposite direction. The
projections of $\Gamma_i^+$ and $\Gamma_i^-$  on the second
coordinate axis are thus strictly increasing in absolute value,
and they do not intersect each other. By choosing $a_1$ through
$a_m$ appropriately we see that we may direct the projections such
that none of the image curves intersect when $t$ is close to zero,
and they all increase in absolute value.  That Condition 1 is
satisfied now follows from Lemma \ref{move}. Condition 2 follows
from Proposition \ref{ext} above, and the proposition is proved.
\end{proof}

\emph{Proof of Theorem \ref{main}:} The result is now an immediate
consequence of Theorem \ref{fb}, Proposition \ref{em}, and the discussion
at the beginning of Section 2. $\square$

\section{Proofs of Theorems \ref{jordan} and \ref{exh}}

To prove Theorem \ref{jordan}, we show that we may change
coordinates such that the conditions in Theorem \ref{main}  are
satisfied.
\begin{lemma}\label{holder}
Let $\mu_1,\mu_2\colon [0,1]\rightarrow\CC$  be
$\mathcal{C}^1$-continuous curves,
 and assume that there is no $t_0\in [0,1]$  such that
$\mu_1(t_0)=\mu_2(t_0)=0$.  Then there exists a set $E\subset\CC$
of measure zero such that for all $c\in\CC\setminus E$  we have
$$
\mu_1(t)+c\cdot\mu_2(t)\neq 0
$$
for all $t\in [0,1]$.
\end{lemma}
\begin{proof}
Let $\{\Omega_j\}$  be a decreasing sequence of open sets in
$[0,1]$  such that
$$
\cap_{j=1}^{\infty}\Omega_j=\{t\in [0,1];\mu_2(t)=0\}.
$$
For each $i\in\NN$  we let
$$
E_i=\{-\frac{\mu_1(t)}{\mu_2(t)};t\in [0,1]\setminus\Omega_i\}.
$$
Then each $E_i$  is a compact curve in $\CC$  and has measure
zero. It follows that the set
$$
E=\cup_{i=1}^\infty E_i
$$
has measure zero.  Assume to get a contradiction that there is a
$c\in\CC\setminus E$  such that
$$
\mu_1(t)+c\cdot\mu_2(t)=0
$$
for a $t\in [0,1]$.  By assumption we have that $\mu_2(t)\neq 0$,
so there is a $\Omega_i$  with $t\notin\Omega_i$.  But then we
have
$$
\mu_1(t)+c\cdot\mu_2(t)=0\Rightarrow
c=-\frac{\mu_1(t)}{\mu_2(t)}\Rightarrow c\in E_i\subset E.
$$
\end{proof}

\begin{lemma}\label{cord1} Let $M$  be a surface as in Theorem \ref{jordan}, and let
$\gamma(t)=(z(t),w(t))$  parametrize $\partial M$.  There exist
coordinates on $\CC^2$  such that
$$
\frac{\partial z}{\partial t}(t)\neq 0,\frac{\partial w}{\partial
t}(t)\neq 0
$$
for all $t$.  Moreover, the new coordinates can be given by a
linear map that is arbitrarily close to the identity.
\end{lemma}
\begin{proof}
Define the following linear map $A\colon\CC^2\rightarrow\CC^2$:
$$
A(z,w)=(z+c_1w,w+c_2z).
$$
Since $\frac{\partial z}{\partial t}(t)$  and $\frac{\partial
w}{\partial t}(t)$  are never zero for the same $t$  it follows
from the previous lemma that the coefficients $c_1$  and $c_2$
can be chosen to prove the lemma.
\end{proof}
\begin{lemma}\label{cord2} Let $M$  be a surface as in Theorem \ref{jordan}.  There exist
coordinates on $\CC^2$  such that there is a point $p\in\partial
M$  with
$$
\pi_1^{-1}(\pi_1(p))\cap\overline M=p.
$$
Also, if we let $z(t)$  parametrize $\pi_1(\partial M)$  near
$\pi_1(p)$  with $z(0)=\pi_1(p)$, we have $\frac{\partial
z}{\partial t}(0)\neq 0$.  Moreover, the new coordinates can be
given by a linear map that is arbitrarily close to the identity.
\end{lemma}
\begin{proof}
Choose coordinates as in the previous lemma.  Since $\overline M$
is compact we have that the restriction of $\pi_1$  achieves its maximum absolute value at a point
$p\in\partial M$.  We may assume that $\pi_1(p)=1$.  Let $P$  be
the set
$$
P=\{p_i\}=\{(1,w_i)\in\partial M\},
$$
which is finite by assumption.  We will show that an arbitrarily
small $\CC$-linear perturbation of $M$ will yield the result. \

First we may perturb the $w$-coordinate such that the set $\{Re
(w_i)\}$  has a unique maximum, say $y_1=Re(w_1)$.  Let
$\gamma(t)=(z(t),w(t))$ parametrize $\partial M$  near $p_1$  with
$\gamma(0)=p_1$.  Make sure that $Re(w'(0))\neq 0$  for the new
$w$-coordinate.  Notice that for small perturbations of the first
coordinate we will still have that $z'(t)\neq 0$.\

Let $(x,y)$ denote the coordinates on the real space
$\RR^2=\{(Re(z),Re(w)\}$  and let $\tilde M$ be the projection of
$M$ onto $\RR^2$.  Let $\tilde\partial$  denote the projection of
$\partial M$  and let $(x(t),y(t))$  be the parametrization of the
projection near $(1,y_1)$.  There exists an $N\in\RR$ such that
$$
\tilde M\subset\{(x,y);|x|<1,|y|<N\}.
$$
Since $y'(0)\neq 0$  and since $x'(t)$  is nonzero for $t\neq 0$
near the origin we have that $\tilde\partial$  near $(1,y_1)$ has
tangent lines with arbitrarily small positive angles $\Theta$ to
the line $l=\{1\}\times\RR$.  And since $\tilde\partial$  is
compact and $y_1$ is the largest value for which $\tilde\partial$
intersects $l$, we may chose $\Theta$  so small that the tangent
line does not intersect $\tilde\partial$ at other points than at a point
$(\tilde x_{1},\tilde y_{1})$  close to $(1,y_{1})$.  Rotate
both axes in $\RR^2$  through the angle $\Theta$  and let the
rotated vectors generate complex coordinates $e_1$  and $e_2$  on
$\CC^2$.  The effect of the rotation is that in the new coordinates the function $Re(z)$
restricted to $\partial$  takes its unique maximum at a point $\tilde p_{1}$.
So we may normalize to get
$Re(\pi_1(\tilde p_1))=1$   and that $Re(\pi_1(q))<1$ for all other
points $q\in\partial M$. Define a function $f\colon
M\rightarrow\CC$  by
$$
f(z,w)=e^{z}.
$$
It follows that $|f(q)|\leq e^1$  for all $q\in\partial M$, with
$|f(q)|=e^1$ only for $q=\tilde p_1$.  It follows then from the maximum
principle that no other point of $\overline M$  projects onto $\pi_{1}(\tilde p_{1})$, and
the proof is finished.
\end{proof}
\begin{remark}
The conclusion of the previous lemma is stable under small linear
perturbations, and this will be used in the proof of Theorem
\ref{torus}. First we locate a point $p$  for one boundary
component, and then we perturb to get another point $q$  for the
other boundary component.
\end{remark}
\emph{Proof of Theorem \ref{jordan}:}  The theorem follows from
the previous lemma and Theorem \ref{main}. \ $\square$

\emph{Proof of Theorem \ref{exh}:} By scaling we may assume that
$M\subset\triangle^2$.  Let $\{K_i\}$  be an exhaustion of $M$ by
polynomially convex compact sets such that $\partial K_i$  is a
set of smooth Jordan curves and $K_i^\circ$  are all diffeomorhic
to $M$.  For a given $K_i$  choose points $p^i_j$, one in each
connected component of $\partial K_i$.  By Theorem 2.3 in
\cite{fr}  there exists an automorphism $\phi\in \Aut(\CC^2)$ such
that the following are satisfied: \

(i) $\phi(p^i_j)=q^i_j\in\partial\triangle^2$  for each $p^i_j$,\

(ii) $|\pi_1(q^i_j)|=1$  for each $q^i_j$, and
$\pi_1(q^i_m)\neq\pi_1(q^i_j)$ for $m\neq j$,\

(ii) $\phi(K_{i-1})\subset\triangle^2$. \

Let $\tilde M$  be the connected component of
$\phi^{-1}(\phi(M)\cap\triangle^2)$  that contains $K_{i-1}$. It
follows from the maximum principle that $\tilde M$  is
homeomorphic to $K_{i-1}^\circ$.  We have that $\tilde M$ contains
a smoothly bounded domain $U$  which is diffeomorphic to
$K_{i-1}^\circ$, such that $K_{i-1}\subset U$, and such that $U$
satisfies the conditions of Theorem \ref{main}. $\square$

\section{Bordered Riemann Surfaces - proof of Theorem \ref{torus}}

In this section we will apply the above results to embed
subsets  $\tilde T$  of the torus with two boundary components properly into $\CC^2$.
To do this, we will use the Weierstrass p-function to embed
$\tilde T$  onto a submanifold of $\CC^2$  satisfying the
conditions in Theorem \ref{main}. \

Let $\omega_1$  and $\omega_2$  be two non-zero complex numbers
that are linearly independent over the real numbers, and let
$X=\{\omega\in\CC;\omega=n\omega_1+m\omega_2, n,m\in\mathbb{Z}\}$.
We obtain a torus $T$  by introducing the following equivalence
relation on $\CC$:
$$
z_1\sim z_2\Leftrightarrow z_1-z_2\in X,
$$
and endowing the quotient space with the obvious complex
structure.  It is known that these are in fact all tori (see for
instance \cite{fk}). Recall that we have the following meromorphic
function defined on the tours (called the Weierstrass p-function):
$$
\varrho(z)=\frac{1}{z^2}+\sum_{\omega\neq
0}(\frac{1}{(z-\omega)^2}-\frac{1}{\omega^2})
$$
We have that $\varrho$  is two-to-one (counted with multiplicity),
and the ramification points are $p_1\sim
0,p_2\sim\frac{\omega_1}{2},p_3\sim\frac{\omega_2}{2},p_4\sim\frac{\omega_1+\omega_2}{2}$.
We claim that for most automorphisms of the torus
$\phi\in\Aut(T)$, the functions $\varrho$  and $\varrho\circ\phi$
separate points. For a complex number $c\in\CC$, let $\phi$  be
the automorphism $\phi(z)=z+c$, and let $z_0\in\CC\setminus X$.
For any $N\in\mathbb{N}$  we let $\Omega_N=\{\omega\in
X;|n|,|m|\leq N\in\NN\}$, and we have that
$$
\sum_{\omega\in\Omega_N}\frac{1}{(z_0-\omega)^2}=\sum_{\omega\in\Omega_N}\frac{1}
{(z_0+\omega)^2}=\sum_{\omega\in\Omega_N}\frac{1}{(-z_0-\omega)^2}.
$$
By letting $N\rightarrow\infty$, we see that
$\varrho(z_0)=\varrho(-z_0)$. We have that $z_0\sim
-z_0\Leftrightarrow 2z_0\in X$, which means that $z_0$ is a
ramification point.  Since $\varrho$
 is two-to-one, we then have that
$$
\varrho(z_0)=\varrho(z_0')\Leftrightarrow z_0'\sim-z_0 \
\mathrm{or} \ z_0'\sim z_0.
$$
Assume that we also have $\varrho(z_0+c)=\varrho(z_0'+c)$.  We get
that
$$
z_0'+c\sim-z_0-c\Leftrightarrow z_0'+2c\sim -z_0\Leftrightarrow
2c\in X.
$$
So let $p,q\in T$  such that $2(p-q)\in\CC\setminus X$. If we let
$f_p(z)=\varrho(z-p)$  and $g_q(z)=\varrho(z-q)$, we get that the
function $\psi\colon T\setminus\{p,q\}\rightarrow\CC^2$ defined by
$$
\psi(z)=(f_p(z),g_q(z))
$$
is one-to-one.  Since $f_p$  and $g_q$  will not share the same
ramification points, we have that $\psi$  is an embedding into
$\widehat\CC^2$. We may sum this up in the following lemma:

\begin{lemma}\label{emb}
Let $p,q\in T$  such that $2(p-q)\in\CC\setminus X$. Then there
exists a proper holomorphic embedding $\psi=(f_p,g_q)\colon
T\setminus\{p,q\}\rightarrow\CC^2$.  Moreover, $f$  and $g$  are
translations of the Weierstrass $p$-function.
\end{lemma}
If we let $D\subset T$  be a disc, we may use this lemma and
Theorem \ref{jordan}  to give a proof of the fact that $T\setminus
D$ can be properly embedded into $\CC^2$.  Let $p,q\in D$ satisfy
the conditions in the lemma, and Theorem \ref{jordan} tells us
that $\psi(T\setminus D)$  can be embedded properly.  In a similar
manner, we will now employ Lemma \ref{emb} to embed the surface in
Theorem \ref{torus} onto a submanifold of $\CC^2$  that satisfies
the conditions of Theorem \ref{main}.

\

\emph{Proof of Theorem \ref{torus}}: We start with the case of two
points.  By Lemma \ref{emb}, for a point $p\in T$ we have that
most choices of $q\in T$ gives us that the map
$\psi=(f_p,g_q):T\setminus\{p,q\}\rightarrow\CC^2$ is a proper
holomorphic embedding.  The only case not covered is the case
where $2(p-q)\in X, p\neq q$.  Let $x_0=p+\frac{q-p}{2}$. Then
$g_{x_0}(p)=g_{x_0}(q)=\a$, and $2(x_0-p)\in\CC\setminus X$. So
the pair $f_p,g_{x_0}$ separate points on $T$, which means that
the pair $f_p,\frac{1}{g_{x_0}-\a}$ separate points.  They do not
share any ramification points, so they furnish the embedding we
are looking for.\

Now to the case of two discs. By choosing the point $p$ close
enough to $\partial D_1$, we may assume that
$$
N=\|f_p\|_{\partial D_1}>\|f_p\|_{\partial D_2}.
$$
Choose a point $x_0\in\partial D_2$ such that $f_p(x_0)\neq 0$.
There exists an $M\in\NN$  such that $|g_{x_0}(z)|<M$  for all
$z\in\partial D_1$, so for $q$  close enough to $x_0$  we have
that $|g_q(z)|<M$  for all $z\in\partial D_1$. If $q$  is close
enough to $x_0$  we have that $|g_q(x_0)|\cdot|f_p(x_0)|>M\cdot
N$.  We get that the function $f_p\cdot g_q$  takes its maximum on
$\partial D_2$.  Let $\tilde T$  denote the embedded bordered
surface $\psi(T\setminus(D_1\cup D_2))$.\

Now, by an argument as in the proof of Lemma \ref{cord2} we may
assume that the conditions in Theorem \ref{main} are satisfied at
a point $p_1\in\psi(\partial D_1)$, and by a change of coordinates
given by a translation in the $z$-coordinate, we may assume that
$\pi_1(p_1)=0$.  Let $F\colon\CC^2\rightarrow\CC^2$ be the map
defined by
$$
F(z,w)=(ze^{h(zw)},we^{-h(zw)}),
$$
where $h\in\mathcal O(\CC)$. Then $F$ is an automorphism of
$\CC^2$, and if $h(0)=0$ it fixes the coordinate axes.  Choose the
function $h$  such that $h(0)=h'(0)=0$, and such that the function
$f_pe^{h(f_pg_q)}$  restricted to $T\setminus(D_1\cup D_2)$  takes
its maximum on $\partial D_2$. Now the projection on the $z$-axis
restricted to the surface $F(\tilde T)$ will take its maximum at a
point $p_2\in F(\psi(\partial D_2))$, and again we may assume that
the condition in Theorem \ref{main} is satisfied at $p_2$.  The
theorem then follows. \

Lastly we consider the case of $D_1$  being a disc, and $D_2$
being a point $q$. We will of course choose a point $p\in D_1$,
and let $\psi(x)=(f_p(x),g_q(x))$  be our first embedding into
$\CC^2$. In this case however, we cannot change coordinates such
that we get a point on the boundary where the projection on the
first coordinate takes its maximum without having any other points
in the fiber. The result will follow from Theorem \ref{main} if we
choose the point $p$
 such that $f_p$  takes its maximum at a point $x_0\in\partial
D_1$, such that $f_p(x_0)$  has no other pre-image. \

Let $z_0\in\partial D_1$  such that $\partial D_1$  is convex near
$z_0$, so that $\mathrm{Im}z$  has a local minimum at $z_0$. We
will show that for small $\d>0$ we can let $p_\d=z_0+i\d$, and
$f_{p_\d}$ will work. \

Let $r_0>0$ such that we may write the Weierstrass p-function as
$$
\varrho(z)=\frac{1}{z^2}+f(z),
$$
for an $f\in\mathcal O(\overline\triangle_{r_0})$. Let
$M=\|f\|_{\overline{\triangle_{r_0}}}$.  We have that
$$
|\varrho(z)|\leq\frac{1}{r^2}+M=\frac{1+Mr^2}{r^2}
$$
for all $|z|=r$  when $r<r_0$.  And for any $k\in\NN$  we have
that
$$
|\varrho(z)|\geq\frac{k^2}{r^2}-M=\frac{k^2-r^2M}{r^2}
$$
for all $z\in\triangle_{r/k}$.  So if $k$  is chosen big enough,
for all $r<r_0$  we have that
$$
(*) \ z\in\triangle_{r/k}\Rightarrow
|\varrho(z)|>\|\varrho\|_{T\setminus\triangle_r}.
$$
With $p_\d$  as above, we have that $\phi_\d(z)=z-p_\d$  is an
automorphism of $T$  such that $\gamma_\d:=\partial(\phi_\d(D_1))$
has negative imaginary part near the point $z_0$  in the new
coordinates. If $\gamma_\d$ passes through $\triangle_{r/k}$, by
$(*)$ we have that to study where $f_{p_\d}$ takes its maximum on
$\partial D_1$, we only have to look at where $\varrho$ takes its
maximum on $\gamma_\d\cap\triangle_r$. Near the origin we may
write $\varrho(z)=\frac{1}{z^2}\cdot g(z)$, where $g\in\mathcal
O(\triangle_r)$, and $g(0)=1$.  So we may let $h^2=g$ and write
$\varrho(z)=(\frac{1}{z}\cdot h)^2$.  It suffices to look at where
$\frac{1}{z}\cdot h$  takes its maximum, or equivalently where
$\varphi(z)=z\cdot\frac{1}{h}=z+O(|z|^2)$ takes its minimum.
Define the set
$$
U_r=\{z\in\triangle_r;\mathrm{Im}z\leq-\frac{r}{2k}\}.
$$
There is a $C\in\RR^+$  (independent of $r$)  such that
$|\varphi(z)-z|\leq C|z|^2$  for all $z\in\triangle_r$, so if $r$
is chosen small enough we have that $\varphi(U_r)$  is contained
in the lower half-plane. \

To show how to choose $\d$, we will now express $\d$  in terms of
$r$, and we let $\d(r)=\frac{r}{k}$.  Near the origin, the curve
$\gamma_{\d(r)}$  can be expressed as a smooth graph over the
x-axis:
$$
y=s(x)-\d(r),
$$
where $s(x)$  is independent of $r$, and $s(0)=s'(0)=0$.  If $r$
is small enough then $s'(x)<\frac{1}{2k}$  for all $x\in (-r,r)$,
so we have that $s(x)<\frac{r}{2k}$  for all $x\in (-r,r)$.  In
other words: $\gamma_{\d(r)}\cap\triangle_r\subset U_r$  if $r$ is
small enough.  Now, $\varphi$ restricted to
$\gamma_{\d(r)}\cap\triangle_r$  has to take its minimum at a
point $z_0$, which means that $\frac{1}{\varphi}$  restricted to
$\gamma_{\d(r)}\cap\triangle_r$  takes its maximum at $z_0$.  We
have that $\varrho=(\frac{1}{\varphi})^2$  (near the origin), and
the square function is 1-1 on $U_r$, so it follows that $\varrho$
(restricted to the boundary) takes its maximum at $z_0$ and that
no other point on the boundary gets mapped to $\varrho(z_0)$.
Since the only ramification point for $\varrho$  near the origin
is the origin, the proof is finished. $\square$ \

\bibliography{biblio}

\end{document}